\documentclass[conference,letterpaper]{IEEEtran}
\IEEEoverridecommandlockouts                        

\usepackage[top=1in, bottom=0.75in, left=0.75in, right=0.75in]{geometry}
\usepackage{graphicx}
\usepackage{amsmath}
\usepackage{amssymb}
\usepackage{eurosym}
\usepackage{booktabs}
\usepackage{xcolor}
\usepackage[]{acronym}
\usepackage[switch]{lineno}
\usepackage[noadjust]{cite} 

\usepackage{soul}
\usepackage{hyperref}

\usepackage[absolute,overlay]{textpos}

\begin{document}

\title{\LARGE \bf Optimal Control of Industrial Assembly Lines*}

\author{Francesco~Liberati, Andrea~Tortorelli, Cesar Mazquiaran, Muhammad Imran and Martina~Panfili
\thanks{*This work has been carried out in the framework of the SESAME project, which has received funding from the European Union’s Horizon 2020 research and innovation programme under grant agreement No 821875. The content of this paper reflects only the author’s view; the EU Commission/Agency is not responsible for any use that may be made of the information it contains.}
\thanks{$^{1}$ The authors are with the DIAG department of Sapienza University of Rome, Via Ariosto 25, 00185, Rome, Italy. Corresponding author: Francesco LIBERATI. Authors' email: 
{\tt\small [surname]@diag.uniroma1.it}.}%
}

\begin{textblock*}{18cm}(1cm,1cm) 
   © 2020 IEEE.  Personal use of this material is permitted.  Permission from IEEE must be obtained for all other uses, in any current or future media, including reprinting/republishing this material for advertising or promotional purposes, creating new collective works, for resale or redistribution to servers or lists, or reuse of any copyrighted component of this work in other works.
\end{textblock*}

\acrodef{ID}{identification}
\acrodef{IP}{integer programming}
\acrodef{MPC}{model predictive control}
\acrodef{MILP}{mixed integer linear programming}

\maketitle
\thispagestyle{empty}
\pagestyle{empty}

\begin{abstract}
This paper discusses the problem of assembly line control and introduces an optimal control formulation that can be used to improve the performance of the assembly line, in terms of cycle time minimization, resources' utilization, etc. A deterministic formulation of the problem is introduced, based on mixed-integer linear programming. A simple numerical simulation provides a first proof of the proposed concept.
\end{abstract}

\begin{IEEEkeywords} assembly line control, industry 4.0, manufacturing, model predictive control.
\end{IEEEkeywords}

\section{Introduction}\label{introduction}

\subsection{Motivation and Objectives}\label{motivation and objectives}
The industrial sector, worldwide, is undergoing a revolution labelled Industry 4.0, whose pillars are digitalisation and automation. The underlying technological developments include: (i) deployment and management of pervasive (wireless) sensor networks, (ii) creation of virtualized environments for simulation, (iii) development of augmented reality applications, (iv) use of big data, (v) full integration of factories with the Information and Communication Technologies (ICT) continuum \cite{2015_davies}. These elements are expected to significantly increase productivity, flexibility, reliability, quality of the manufactured products and revenues \cite{2015_boston}, \cite{2012_fiaschetti}.

The mentioned technological developments can be implemented in all the industrial life cycle phases (i.e. procurement, production, distribution, sales). In this paper, the focus will be on the production phase and, in particular, on assembly lines. The objective is to present an approach for increasing flexibility of industrial assembly lines, through an optimization framework able to, e.g., minimize the cycle time and perform the re-scheduling of line activities, in presence of unexpected events captured by the sensor networks, and possibly elaborated by means of big data analytics techniques. These aspects will be further detailed in the next sections.
%
\subsection{Literature Review}\label{literature review}
The optimization of industrial assembly lines is a well known problem in the literature and many formulations have been proposed. The literature will be reviewed with a focus on the Assembly Line Balancing Problem (ALBP).

An assembly line consists of a set of workstations, each performing different operations (tasks) on items moving on the line according to a maximum or average time, called the cycle time. The ALBP consists in the optimal assignment of tasks to workstations, while always satisfying given constraints and requirements. Several formulations of the ALBP have been proposed.

In \cite{1986_baybars} two main families of ALBP problems have been identified: Single ALBP (SALBP) and Generalized ALBP (GALBP). The SALBP is characterized by the production of a unique element, a fixed common cycle time, a serial one-sided line, precedence constraints and deterministic times among other constraints \cite{2007_boysen}. The GALBP formulation relaxes some of the SALBP assumptions for addressing more realistic scenarios. In particular, GALBPs allow to consider multi-model process, zoning constraints, delays, parallel stations and more complex layouts among others. In this latter class of problems there are the U-shape lines (UALBP) and the Mixed-Model lines (MMALBP).

A further classification can be performed based on the optimization criteria considered, which characterize the so-called problem-type. A Type-1 problem consists in optimizing the number of workstations for a given task and a given cycle time. A Type-2 problem, instead, consists in the minimization of the cycle time for a given number of workstations. In Type-E problems, the objective function is the efficiency of the assembly line. In Type-F problems the objective is to find a feasible balance between the number of workstations and cycle time which are both fixed \cite{2006_becker},  \cite{2018_kamarudin}.
The above standard problem types can be further particularized to specific applications by customizing the mathematical formulation of the optimization problem.

In \cite{2016_razali}, for example, it is presented an MMALBP Type-2 problem characterized by three different product models in which the resources used in the assembly line and the product rate variation are minimized.
In \cite{2017_rabbani}, the authors solve a MMALBP addressing simultaneously the line balancing and sequencing problems, comparing the performances of two algorithms: the multi-objective particle swarm optimization (MOPSO) and the non-dominated sorting genetic algorithm (NSGA-II); results show that the latter has better performances.
In \cite{2018_kamarudin}, the authors address a SALBP-1 problem with resources constraints. Three different objective functions are taken into account: minimization of workstations, machines used and the number of multi-skilled workers.
In \cite{2018_pereira}, the authors show that the quality of SALBP solutions decreases when products or resources increase. This is achieved by considering the production of multiple products in a single line, and the use of special machinery needed to perform some tasks.
In \cite{2015_gyulai}, the authors propose a \ac{MILP} problem to plan the allocation of jobs and workers at workstations, to minimize costs of setup, personnel and early or late delivery. Capacity requirements are estimated via linear regression on historical data, which makes the model more robust.
In \cite{2018_manzini}, the authors propose an integrated approach to optimal design of reconfigurable assembly systems, taking into account future scenarios of demand, and the consequent need of reconfiguring the assembly system. The paper proposes an ``assembly system configuration tool", ``an assembly cell configuration tool", based on a set of logic constraints capturing the cells' dynamics, a ``production planning tool" (based on \ac{MILP}) and a ``reconfiguration planning tool" (again based on \ac{MILP}). The overall objective is to minimize costs. 
In \cite{2017_gyulai} the authors define a two level problem, for first planning the production lots' size (by matching capacity with demand), at a lower level, scheduling the execution of the planned tasks. The two formulations are based on \ac{MILP}, the overall objective is to minimize operator costs and penalties on early and late delivery. Finally, \cite{2020_koskinen} provides an \ac{IP} model and a heuristic for optimizing the planning of a multi-line assembly system for Printed Circuit Board manufacturing. The optimization goal is to minimize the total tasks tarty times (i.e., the accumulated delays) and balance lines. The heuristic is shown to scale well also in real scenarios. 
%
\subsection{Paper Contributions}\label{contribution}
Building on the literature, in this paper we address the problem of optimally planning and controlling the execution of tasks in an assembly line. 
Compared to the literature, the present model is more detailed and flexible, in the sense that it captures all the main relations and constraints characterising resources (workers, tools, components, transportation vehicles, etc.) and tasks. Such constraints are often dynamical in nature, and need to be captured to have a fine grained and applicable task controller. For example, most of the time in literature the transportation of resources across the assembly line is not modelled in detail as well as the qualification/authorization requirements of workers to be assigned to tasks, etc. This aspect is especially relevant in complex operational environments, where given tasks can only be executed by given workers. All these aspects are crucial for an optimal planning and control of operations. The model proposed in this work allows to address all of these aspects allowing to provide human operators with a valuable tool in support to decision making processes. 
\subsection{Outline}\label{outline}
The remainder of this paper is organized as follows. Section II presents the addressed problem and the reference scenario. Section III presents an outline of the mathematical formulation of the task planning and control problem. Section IV presents a simple and illustrative numerical simulation. Section V concludes the work.
\section{Problem Description and Reference Scenario} \label{problem description and reference scenario}
We consider the problem of optimally scheduling  and then controlling in real time the execution of a given set of tasks in a given assembly line. %
The tasks' characteristic (such as duration, needed input resources and resulting output resources) are known inputs of the problem. Also, the layout and the characteristics of the assembly line are known, for example in terms of available workstations, technical specifications of the workstations, tools and personnel available, etc.

The first scenario in which the proposed algorithm could be applied regards the static optimization of a given assembly process, composed of a high number of assembly tasks. It is not uncommon in fact that, in practice, the assembly schedule followed by the operators is a sub-optimal one, especially in the most complex cases when the single human operator or a team of scheduling operators cannot tame the combinatorial complexity of the scheduling problem. This first scenario is thus an optimal planning of the scheduling activities.

A second scenario is the real time control of the planned schedule, to make sure that the progress is in line with the schedule, and to provide real-time re-planning functionality in case of anomalies. In this second scenario thus, the algorithm can be used to provide real time decision support to the operator responsible for the scheduling activities.

\begin{figure}[tb]
\centering
\includegraphics[width=1\columnwidth]{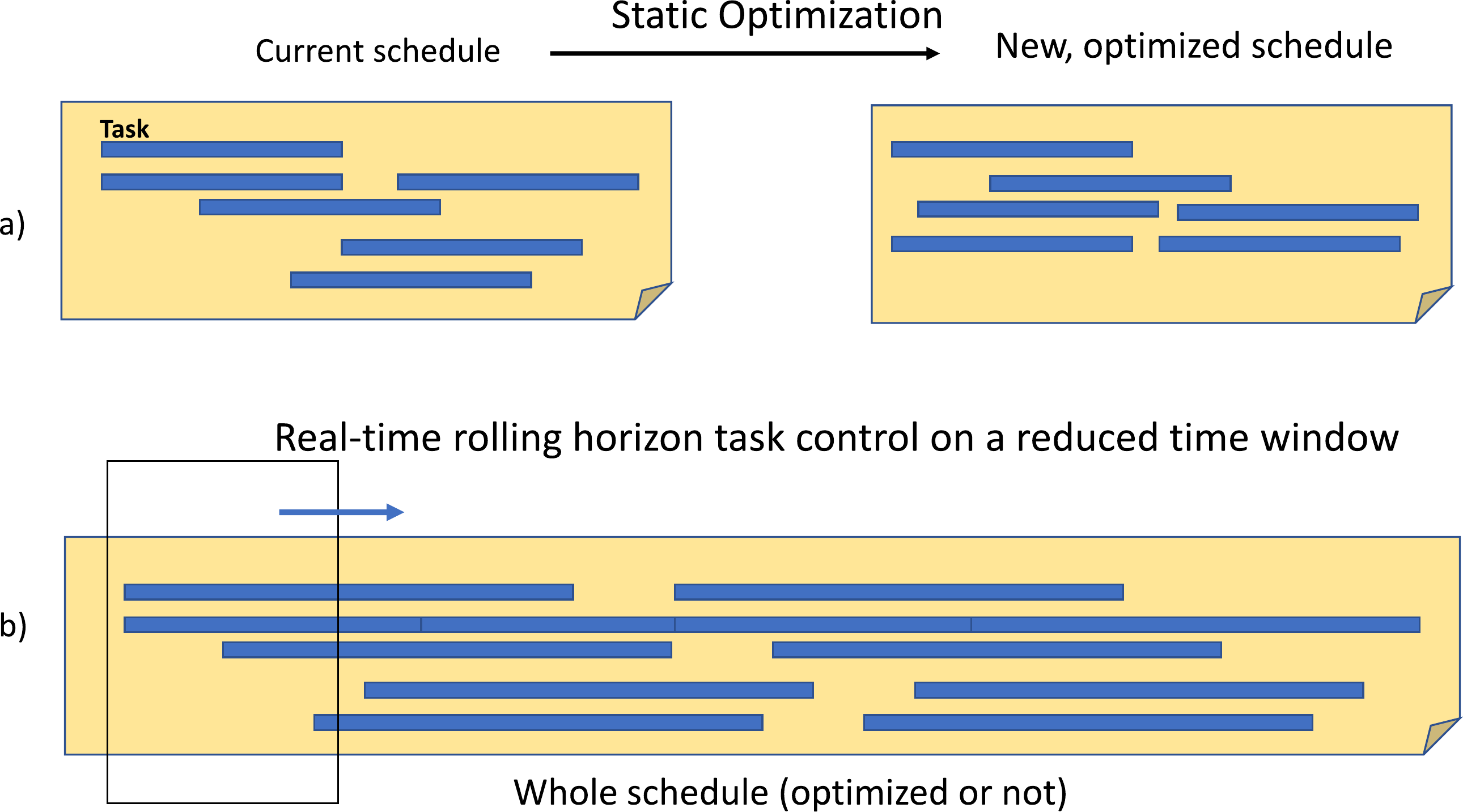}
\caption{Two usage scenarios: a) full plan optimization; b) real time planning control/re-optimization.}
\label{f1}
\end{figure}

The two scenarios are depicted in Fig. \ref{f1}. The first scenario (Fig. \ref{f1}.a) aims at optimizing the whole task planning. Optimization is done once (i.e., static optimization). In case of very complex schedules, the problem could require powerful hardware or could be even intractable (for the high number of variables in the model). In the second scenario (Fig. \ref{f1}.b), once an overall planning of the tasks execution is given (either optimal or not), the proposed algorithm runs in real time to enforce the scheduling, and to correct it in case that anomalies happen, which requires the real-time computation of alternative solutions. In complex scenarios, the real time optimization would be performed on a reduced time window, and not on the overall time span of the whole planned schedule (which again could be computationally intractable).

In this paper, we outline the mathematical model underlying both formulations. The simulations will refer to a simple illustrative case, in which the point of view of the first usage scenario can be directly adopted.
%
\section{Problem Formulation} \label{problem formulation}
In this section we first model the objects present in the problem (tasks, workstations, transportation vehicles, etc.), and then outline the mathematical formulation of the proposed optimal task control problem.
For reason of space, we limit in the following to a high level presentation of the mathematical framework, and we omit to detail some of the constraints (which would require heavy mathematical notation to be precisely defined). The full, detailed formulation will be subject of a future publication.
\subsection{Model of Tasks, Resources and Workstations}
\subsubsection{Tasks} $\mathcal{T}$ is the set of all tasks to be planned/controlled. Tasks are characterised by a duration $d_t$, an earliest allowed starting time, $S_t$, and a latest allowed finish time, $F_t$. 
\subsubsection{Resources} $\mathcal{R}$ is the set of all the resources (tools, materials, people, etc.) available in the assembly line. Let $\mathcal{R}^t$ denote the set of resources types (two or more resources can be of the same type). Example of resource types are: $\mathcal{C}$ consumables, i.e., resources which are consumed by the tasks (e.g., water, power, chemicals, parts to be assembled, etc.), $\mathcal{I}$, instruments/tools/operators needed to execute a task (freed after the task ends), $\mathcal{V}$, the set of transportation vehicles, etc.
Tasks consume and produce resources. We denote with $R^{i}_{t}$ the input resources required by task $t$, and with $R^{o}_{t}$ the output resources/products produced by the task $t$.
\subsubsection{Workstations} $\mathcal{W}$ is the set of workstations, i.e., the specific machines and/or sectors of the assembly line where the tasks are worked. $\mathcal{W}_t$ denotes the set of workstations that can work task $t\in\mathcal{T}$. We model the workstation as a system which is composed by a processing unit, which has an internal capacity (which limits the maximum amount of tasks that can be worked in  parallel), and a buffer for storing the input and the output resources.
A particular type of workstation are the source and storage points of material where the resources are stored.
\subsubsection{Factory Layout}
We model the overall factory as a graph $\mathcal{G}=\{\mathcal{W},\mathcal{E}\}$, where the set of the graph's nodes $\mathcal{W}$ is the set of workstations and $\mathcal{E}$ is the set of transportation paths connecting the workstations. By convention, we let $(w_1,w_2)\in\mathcal{E}$, only if $w_1<w_2$.
\subsection{Definition of the State Variables}
The state of the assembly line at a given time $k$ is given by the collection of the following variables.
\begin{enumerate}
    \item \textit{Workstation state variables}. 
    \begin{enumerate}
        \item $o_{w,k}$ the \textit{occupancy level} (i.e., number and type of tasks in execution) of workstation $w$ at time $k$. The occupancy level dictates the number of tasks that can be accepted at each given time by the workstation;
        \item $R_{w,k,r}$, the level of the input/output buffer of resources at workstation $w$ at time $k$, for resource type $r$. This dictates the number/type of tasks that can be started at the workstations (since tasks need input resources to execute).
    \end{enumerate}
    \item \textit{Resources state variables}.
        \begin{enumerate}
            \item \textit{Resources Location}. Variable $l_{r,k,(w_1,w_2)}$ denotes the location of resources in the factory (i.e., over the graph $\mathcal{G}$). It is $l_{r,k,(w,w)}\in\{0,1\}$ (i.e., \textit{Boolean} variable), $\forall w\in\mathcal{W}$, with $l_{r,k,(w,w)}=1$ if and only if the resource is at the workstation $w$. It is instead $l_{r,k,(w_1,w_2)}\in[0,1]$ (i.e., \textit{continuous} variable) $\forall (w_1,w_2)\in\mathcal{E}$, with $l_{r,k,(w_1,w_2)}=l>0$ meaning that the resource is on path $(w_1,w_2)$, with $l$ the percentage of the path travelled.
        \end{enumerate}
    \item \textit{Task state variables}.
        $\{i_{t,k},e_{t,w,k},f_{t,k}\}$, where $i_{t,k},e_{t,w,k},f_{t,k}=1$ if and only if, respectively, task $t$ is \textit{idle}, \textit{executing} (at $w$) or \textit{finished} at time $k$.
\end{enumerate}
%
\subsection{Definition of the Decision Variables}
The control variables are:
\begin{enumerate}
    \item $a_{t,w,k}$ (``assign"), a Boolean variable equal to one if and only if task $t$ is assigned to workstation $w$ at time $k$. A second task control variable is $s_{t,k}$, which is equal to one if and only if task $t$ is started at time $k$. It is:
    \begin{equation}\label{transition_r}
        s_{t,k}=\sum_{w\in\mathcal{W}_t}a_{t,w,k}, \ \forall t, k.
    \end{equation}
    \item $p_{r,v,k}$ (``place"), a Boolean variable equal to one if and only if resource $r\in \mathcal{R}$ (e.g., a driver, a tool, a semi-finished component) is placed on vehicle $v$ at time $k$.
    \item $\overrightarrow{g}_{v,k,w_1,w_2}$ (``goto"), a Boolean variable equal to one if and only if vehicle $v$ moves from $w_1$ to $w_2$ at $k$;
    \item $\overleftarrow{g}_{v,k,w_1,w_2}$ (``goto"), a Boolean variable equal to one if and only if vehicle $v$ moves from $w_2$ to $w_1$ at $k$.
\end{enumerate}
%
\subsection{Definition of the State Dynamics}
We specify in the following the laws according to which the state of the system evolves in time.
\begin{enumerate}
    \item \textit{Workstations' occupancy level}. The following equation captures the dynamics of the occupancy level at a workstation at a given time:
    \begin{equation}
        o_{w,k} = \sum_{t\in\mathcal{T}}c_{t}e_{t,w,k}, \ \forall w, k.
    \end{equation}
    The occupancy level can capture different aspects depending on the task and the workstation types. For example, it can be simply the number of tasks executing at the workstation (in which case, it is $c_{t}=1$), or the volume occupied, or the weight (in which case $c_{t}$ is a volume or weight, etc.).
    \item \textit{Dynamics of the workstation input/output inventories}.
    The level of resources at a workstation is equal to the level of resources currently located at the workstation, minus the resources consumed by starting tasks, plus the resources generated/freed by the finishing tasks;
    \item \textit{Dynamics of the vehicles' locations}: Transportation vehicles move on the graph of the workstations, to transport input/output resources across the workstations, to feed the tasks. The general equation of their dynamics is, for all $v,k,(w_1,w_2)\in\mathcal{E}$
    \begin{equation}\label{lkp1}
    \begin{split}
        &(\overrightarrow{g}_{v,k,w_1,w_2} + \overleftarrow{g}_{v,k,w_1,w_2})(l_{v,k+1,(w_1,w_2)} +\\ & - l_{v,k,(w_1,w_2)}) = Tv_{v}(\overrightarrow{g}_{v,k,w_1,w_2} - \overleftarrow{g}_{v,k,w_1,w_2}),
    \end{split}
    \end{equation}
    where $T$ is the sampling time, $v_v$ the velocity of the vehicle. The equation is effective only on the links $(w_1,w_2)$ for which $\overrightarrow{g}_{v,k,w_1,w_2} + \overleftarrow{g}_{v,k,w_1,w_2}=1$, where it is $l_{v,k+1,(w_1,w_2)} = l_{v,k,(w_1,w_2)} +  Tv_{v}(\overrightarrow{g}_{v,k,w_1,w_2} - \overleftarrow{g}_{v,k,w_1,w_2})$. This is a non-linear (quadratic) equation, which can be exactly linearized (as shown in \cite{2013_sherali}) to have faster solving times.
    \item \textit{Dynamics of the task state}. 
    Task dynamcis is specified by a set of equations which dictates the transition from one state to the other, based on the decision variables. For example, the state transits to ``execution" when $s_{t,k}=1$, i.e., when $a_{t,w,k}=1$ for some $w$:
    \begin{equation}
        a_{t,w,k}\leq e_{t,w,k}, \ \forall t, w, k.
    \end{equation}
    Then, let $d_{t,k}$ denote the number of time steps left to complete task $t$ at time $k$. 
    If a task is in the execution state, i.e., $e_{t,w,k}=1$ for some $w$, then it shall remain in the execution state as long as $d_{t,k}>0$.
    If $d_{t,k}=0$, then the task shall no longer be in the execution state, i.e.,
    \begin{equation}
        e_{t,w,k}\leq d_{t,k}, \ \forall t, w, k.
    \end{equation}
    and will go to the ``finished" state.
    Finally, it is not possible to go back from the finished state:
    \begin{equation}
        f_{t,k}\leq f_{t,k+1}, \ \forall t, k.
    \end{equation}
\end{enumerate}
%
\subsection{Other Problem Constraints}
The remaining problem constraints are as follows:
\begin{enumerate}
    \item \textit{Unique task starting time}. Only one starting time can be chosen for every task, between the earliest allowed starting time and the latest feasible one
    \begin{equation}
        \sum_{k=S_t}^{F_t-d_t}s_{t,k}=1, \ \forall t.
    \end{equation}
    \item \textit{Unique workstation of execution}. A task can be in execution at most at one workstation
    \begin{equation}
        \sum_{w\in\mathcal{W}}e_{t,w,k}\leq 1, \ \forall t, k.
    \end{equation}
     \item \textit{Unique task assignment to workstation}. A task can be assigned to one and only one workstation, only once.
    \begin{equation}
        \sum_{k=S_t}^{F_t-d_t}\sum_{w\in\mathcal{W}}a_{t,w,k}=1, \ \forall t.
    \end{equation}
    \item \textit{Task time dependencies constraints}. Tasks can be in one of the following temporal dependence relations:
    \begin{enumerate}
        \item \textit{Finish to Start}: every task $t'$, predecessor of task $t$, must finish before $t$ can start:
        \begin{equation}
            d_{t'} + \sum_{k=S_{t'}}^{F_{t'}-d_{t'}}ks_{{t'},k} \leq  \sum_{k=S_t}^{F_t-d_t}ks_{t,k}.
        \end{equation}
        Note that the left hand side represents the finish time of task $t'$, while the right hand side is the start time of task $t$;
	    \item \textit{Start to Start}: every task $t'$, predecessor of task $t$, must start before $t$ can start:
        \begin{equation}
            \sum_{k=S_{t'}}^{F_{t'}-d_{t'}}ks_{{t'},k}\leq \sum_{k=S_t}^{F_t-d_t}ks_{t,k};
        \end{equation}
	    \item \textit{Finish to Finish}: every task $t'$, predecessor of task $t$, must finish before $t$ can finish:
        \begin{equation}
            d_{t'} + \sum_{k=S_{t'}}^{F_{t'}-d_{t'}}ks_{{t'},k} \leq d_{t} + \sum_{k=S_t}^{F_t-d_t}ks_{t,k};
        \end{equation}
	    \item \textit{Start to Finish}: every task $t'$, predecessor of task $t$, must start before $t$ can finish:
        \begin{equation}
            \sum_{k=S_{t'}}^{F_{t'}-d_{t'}}ks_{{t'},k} \leq d_{t} + \sum_{k=S_t}^{F_t-d_t}ks_{t,k}.
        \end{equation}
    \end{enumerate}
    \item \textit{Task state constraints}. A task can be in one and only one state at every given time.
    \begin{equation}\label{one_state}
        i_{t,k} + \sum_{w\in\mathcal{W}_t}e_{t,w,k} + f_{t,k} = 1, \ \forall t, k. 
    \end{equation}
    \item \textit{Workstations occupancy constraints}.
    \begin{equation}
        o_{w}^{min}\leq o_{w,k}\leq o_{w}^{max}, \ \forall w, k.
    \end{equation}
    \item \textit{Workstation resource inventory constraints}.
    \begin{equation}
        R^{min}_{w,r}\leq R_{w,k,r}\leq R^{max}_{w,r}, \ \forall w, k, r\in\mathcal{R}^t.
    \end{equation}
    The above constraints also make sure that a task is assigned to a workstation only if there are always enough input resources and enough ``space" to store the output resources.
    \item \textit{Vehicles localization constraints}.
    A complex set of constraints is included to enforce the proper movement of vehicles along the paths connecting the workstations. 
    \item \textit{Resources localization constraints}.
    Another set of constraints is needed to model the movement of resources across the assembly line. In short, resources cannot move if they are not associated to a vehicle. On the contrary, when they are associated to a vehicle, their location variables will coincide with the location variables of the vehicle they are associated with.
\end{enumerate}
%
\subsection{Objective Function}
The objective function is written as the convex combination of several terms, reflecting different optimization goals:
\begin{equation}
    V=\sum_{i=1}\alpha_iV_i,
\end{equation}
with $\alpha_i\geq 0$ and $\sum_{i}\alpha_i=1$. 

Possible relevant terms are discussed in the following.
\begin{enumerate}
    \item $V_1$ (\textit{minimization of task completion time}). We pose:
    \begin{equation}
        V_1=\tau
    \end{equation}
   where $\tau$ is an auxiliary variable greater or equal than all the finish times.
    \begin{equation}
        \tau\geq d_{t} + \sum_{k=S_t}^{F_t-d_t}ks_{t,k} \ \ \forall t
    \end{equation}
    Hence, when minimized, $\tau$ represent the latest task finish time, and the inclusion of $V_1$ minimizes it;
    \item $V_2$ (\textit{minimization and balancing of the input/output resource inventory}). %
    \begin{equation}
        V_2 = \sum_{w,k,r}\big(R_{w,k,r} - R_{w,r}^{ref}\big)^2
    \end{equation}
   where $R_{w,k}^{ref}$ is the reference level of the resources inventories (if put to zero, inventory is minimized);
   \item Cost minimization. Cost minimization can be easily integrated, for example to take into account the operational cost of utilizing workstations, the cost of storing inventory, the penalty cost of early or late delivery, etc.
\end{enumerate}
In the following simulations we will consider only the minimization of the overall tasks completion time (i.e., the cycle time).
\section{Simulations}
\subsection{Definition of the Simulation Scenario}
We present a simple simulation to provide a proof of the concept. 
The analysis of a more complex scenario is demanded to a future publication. 
Simulations are performed with the open-source technical computing language Julia, version 1.3.1 (https://julialang.org/). The optimization problem is written using the Julia JuMP modelling package \cite{dunning2017jump} and solved using Gurobi optimizer \cite{gurobi}, 
on an Intel I7, 8GB RAM machine running Windows 10. 

\begin{figure}[tb]
\centering
\includegraphics[width=0.8\columnwidth]{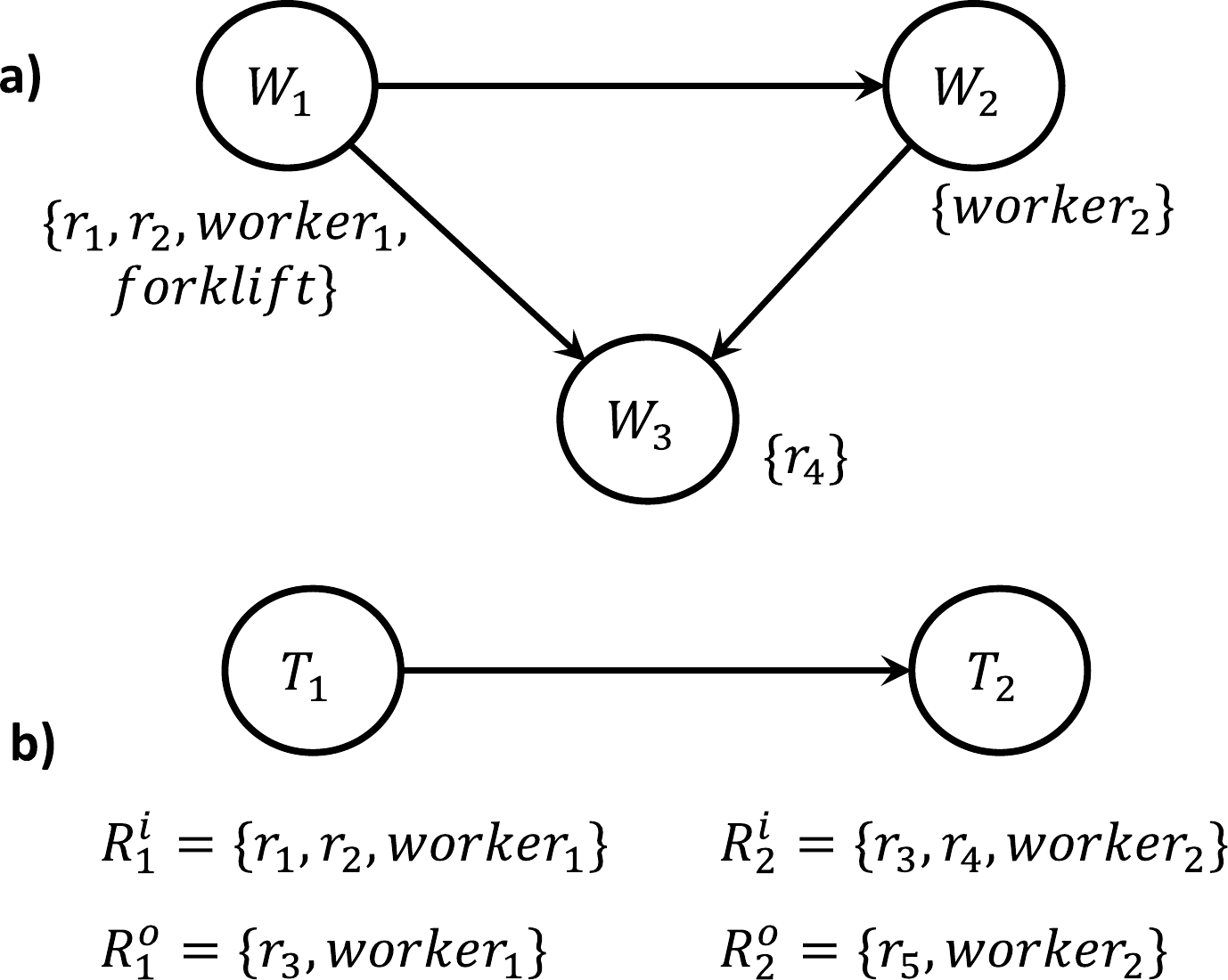}
\caption{Workstations (a) and tasks (b).}
\label{f2}
\end{figure}

The simulation scenario is as follows. We consider an assembly line with three workstations and a job of two tasks, as illustrated in Fig. \ref{f2}.
Figure \ref{f2}.a reports the topology of the paths linking the workstations, and the initial state of the workstations' resource buffers (i.e., the distribution of the resources in the plant). Figure \ref{f2}.b reports the task dependency graph and the inputs needed to execute each task, as well as the outputs produced by the tasks.

The simulation is on a time horizon of $20$ time steps. The tasks last $3$ time steps, can start at any time and at any workstation (provided that the input resources are available at that time/workstation). One vehicle is present in the scenario. A constant velocity is assumed for it.

\subsection{Simulation Result}
\begin{figure}[tb]
\centering
\includegraphics[width=1\columnwidth]{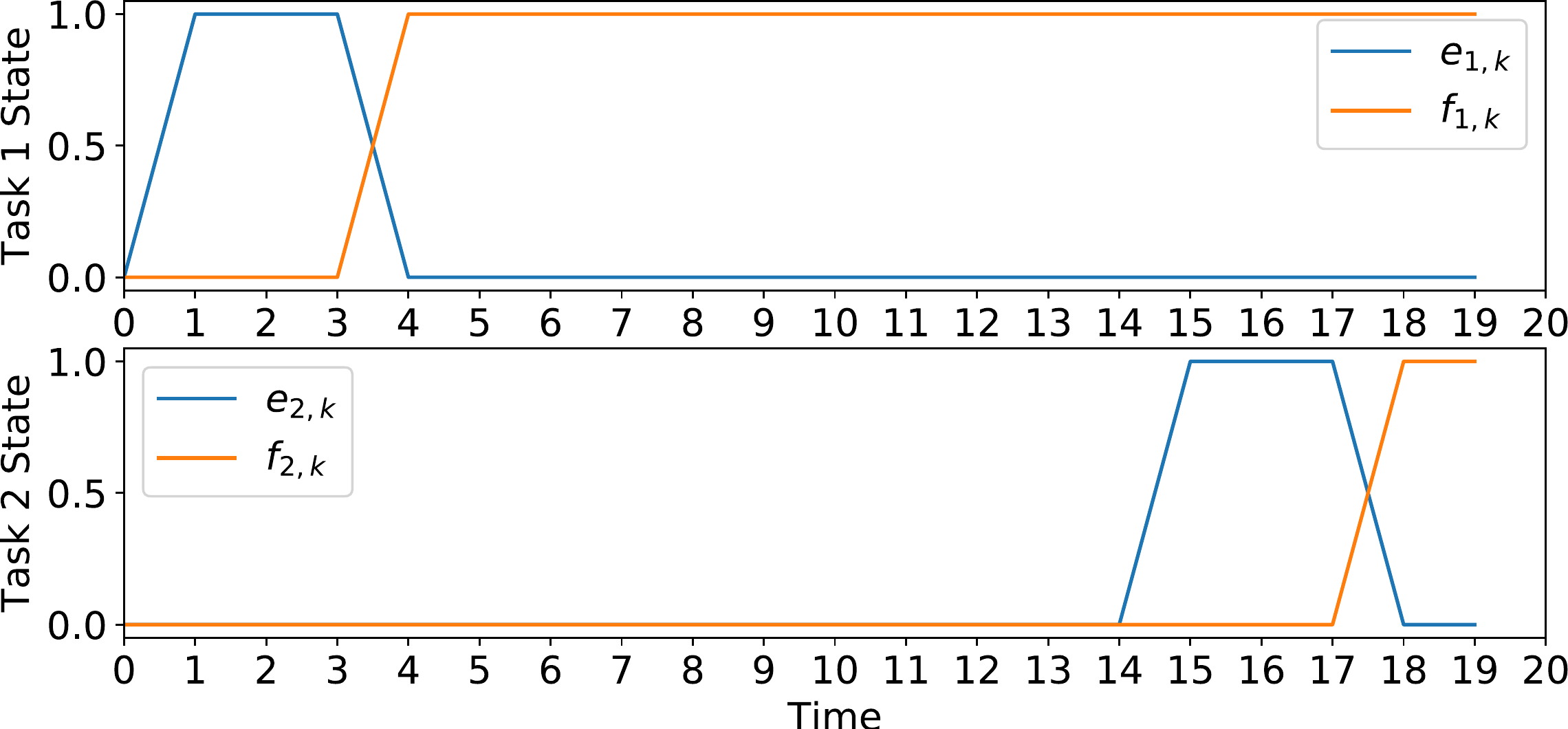}
\caption{Tasks' state.}
\label{f3}
\end{figure}
Figure \ref{f3} reports the tasks' state: task 1 is immediately executed, since it has available at $w_1$ all the needed inputs. Task 2 instead needs $r_3$ (which is generated at $w_1$ as output of task 1), worker 2, who is at $w_2$, and $r_4$, which is at $w_3$. 
\begin{figure}[tb]
\centering
\includegraphics[width=1\columnwidth]{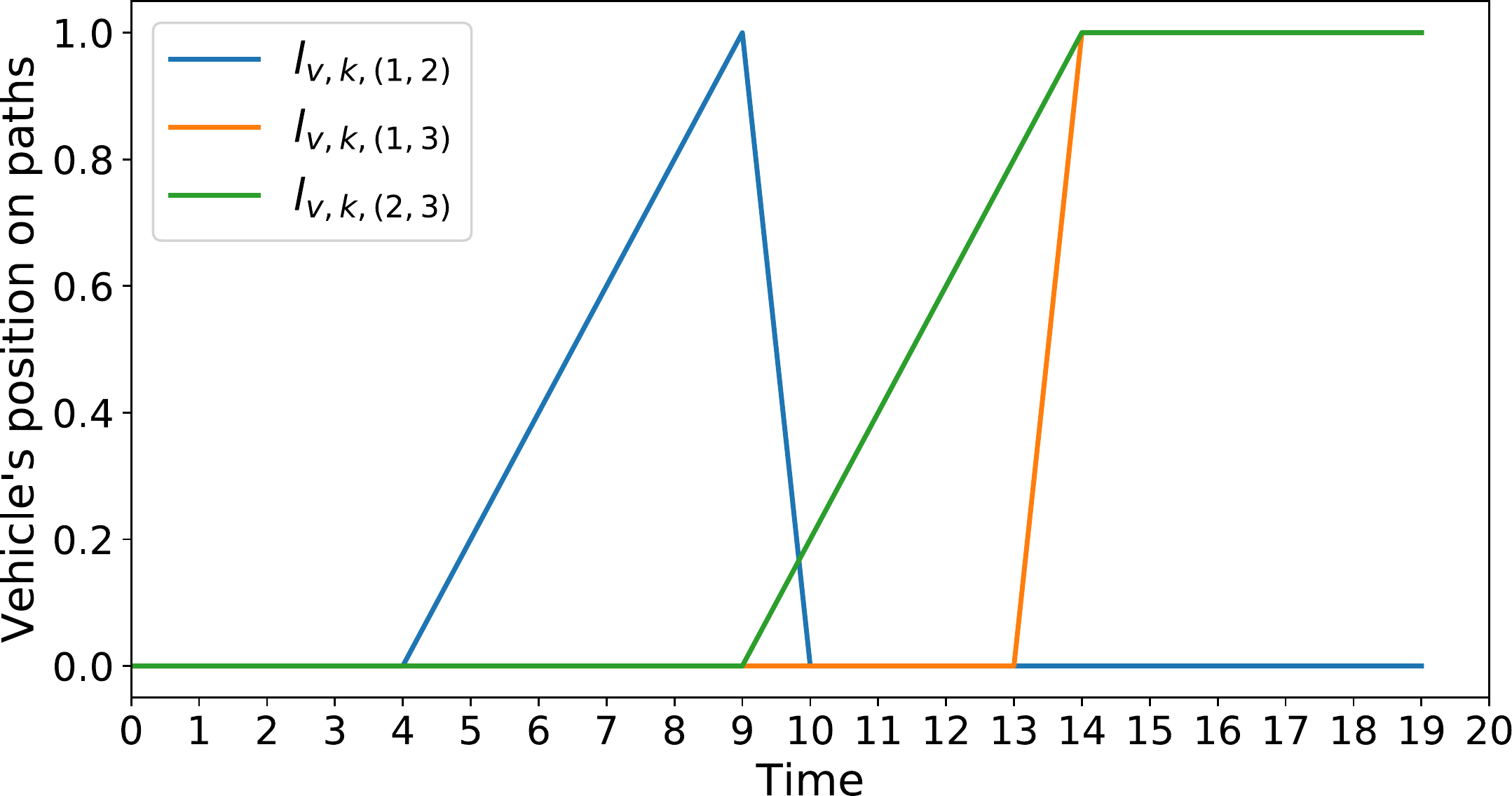}
\caption{Vehicle's movements.}
\label{f4}
\end{figure}
Figure \ref{f4} reports the movements of the vehicle (i.e., its location along the paths connecting the workstations). After $r_3$ is generated, the vehicle carries it from $w_1$ to $w_2$ (i.e., $l_{v,k,(1,2)}$ grows from $0$ to $1$). Once at $w_2$, the vehicle carries also worker 2 and brings both $r_3$ and worker 2 at $w_3$, where then all the resources are available to start task 2. Notice that, when the vehicle is at a workstation, it is also at the end of all the paths incoming into that workstation (e.g., when the vehicle is at $w_3$, both $l_{v,k,(1,3)}$ and $l_{v,k,(2,3)}$ are equal to one).
\begin{figure}[tb]
\centering
\includegraphics[width=1\columnwidth]{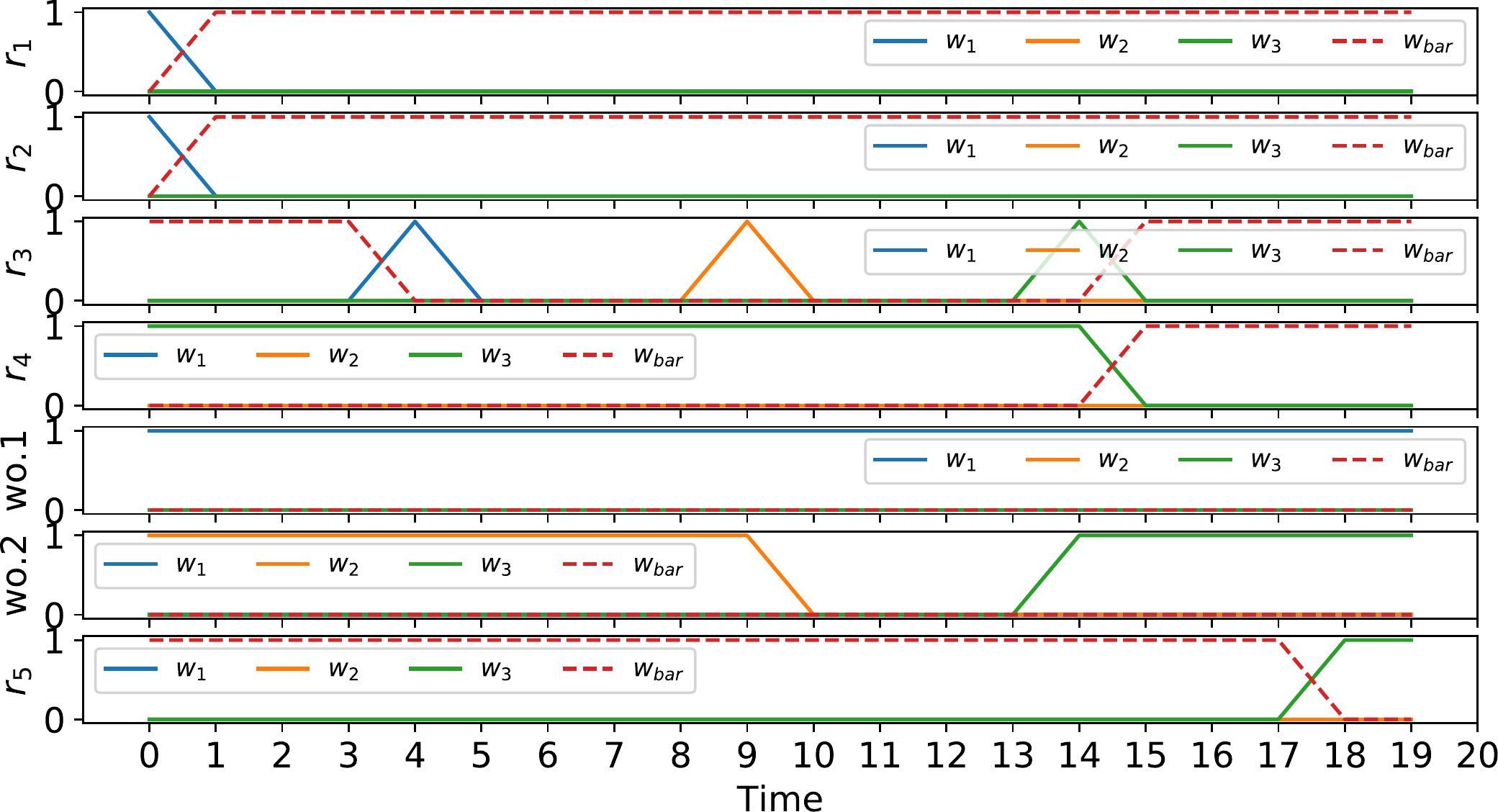}
\caption{Resources' movements.}
\label{f5}
\end{figure}
Figure \ref{f5} reports the location of resources in time ($w_{bar}$ is a fictitious workstation that stores the resources still to be created and the consumed ones). It is seen that the location of resources is consistent with that of the vehicle and with task execution. Notice for example that $r_3$ is created at time $4$ at $w_1$, then at time $9$ it reaches $w_2$, then at time $14$ it finally arrives at $w_3$, where it is consumed by task 2.
The solution found was the one minimizing the job finish time.
Simulation time was $10.75$ seconds. There are $3043$ variables in the model.

\section{Conclusions}
\label{conclusions}
This paper has presented an optimal control-based approach for the scheduling of tasks in an assembly line. The derived optimization problem is based on mixed-integer linear programming and it allows a fine-grained control and optimization of resources (including, personnel, vehicles, input/output resources, etc.) and tasks' execution time. A simple proof of concept simulation has been presented on an illustrative scenario in which the minimization of the total job time is sought. Though there is an inherent combinatorial complexity arising from the presence of binary variables, the linear formulation derived allows the model to be scaled on realistic scenarios. 
Future works will regard the extensive simulation of the proposed framework on realistic scenarios, also using state of the art modelling and simulation tools. The authors are currently investigating model predictive control and deep learning approaches for addressing the computational complexity that may arise in realistic scenarios. Indeed, in some scenarios the sub-optimality of solutions may be acceptable given complexity and real-time requirements. Finally, while the present paper has mostly shown the application of the proposed optimization framework to task planning, the focus of future works will also be more on real time task control (i.e., rescheduling) in presence of uncertainties and disturbances, for which model predictive approach appears a promising candidate methodology.

\section*{Acknowledgment}
The authors thank the CRAT team for the valuable contribution to the design of the algorithm, and all the members of the SESAME H2020 European project for the fruitfull discussions.
\bibliographystyle{myIEEEtran.bst}
\bibliography{bibliografia}
\end{document}